\documentclass[a4paper,12pt,twoside]{article}
\usepackage{amssymb}
\usepackage[final]{lpretex}
\usepackage{title}
\usepackage{a4}
\usepackage{amscd}
\usepackage{amsmath}
\usepackage{bbm}

\input xy
\xyoption{all}

\EndOfDump

\def\DHpartformat#1{(#1) }

\LBsetstyleof{partno}{arabic}

\DocVersion[Changed to LaTeX]{1}

\let\define\def

  \def\C {{\mathbb C}}

\def\GG {{\mathbb G}}   
  
  \def\P {{\mathbb P}}

\def\Z {{\mathbb Z}} 

\define \n {\mathbb N}
\define \z {\mathbb Z}
\define \q {\mathbb Q}
\define \PP {\mathbb P}

\def\sA {{\Cal A}} \def\sB {{\Cal B}} \def\sC {{\Cal C}}
 \def\sE {{\Cal E}} \def\sF {{\Cal F}}

 \def\sN {{\Cal N}} 
 
  \def\sU {{\Cal U}}

\define \cN {\Cal N}
\define \cf {\Cal F}
\define \cg {\Cal G}
\define \cE {\Cal E}
\define \ce {\Cal E}
\define \cc {\Cal C}
\define \cV {\Cal V}
\define \cA {\Cal A}
\define \cK {\Cal K}
\define \cO {\Cal O}
\define \cF {\Cal F}
\define \cn {\Cal N}
\define \cI {\Cal I}
\define \sP {\Cal P}
\def\tA {\widetilde{\Cal A}}

\define \x {\xi}
\define \y {\eta}
\define \G {\Gamma}
\define \r {\rho}
\define \w {\omega}

\def \trho {\widetilde {\rho}}

\def \tp {\widetilde{\mathbb P}}
\define \tH {\widetilde H}
\define \tG {\widetilde{\Gamma}}
\define \tW {\widetilde W}
\define \tF {\widetilde F}
\define \tm {\widetilde m}
\define \St {\widetilde S}
\define \Xt {\widetilde X}
\define \tS {\widetilde S}
\define \tpsi {\widetilde \psi}
\define \tL {\widetilde L}
\define \tE {\widetilde E}
\define \tl {\widetilde l}
\define \tA {\widetilde A}
\define \tom {\widetilde\omega}
\define \tT {\widetilde T}
\define \tB {\widetilde B}
\define \tf {\widetilde f}
\define \tsA {\widetilde{\sA}}

\define \tM {\widetilde M}
\define \tphi {\widetilde{\phi}}
\define \trho {\widetilde{\rho}}
\define \tR {\widetilde R}
\define \tp {\widetilde p}
\define \tq {\widetilde q}
\define \tc {\widetilde c}
\define \tsF {\widetilde {\sF}}
\define \tsN {\widetilde {\sN}}
\define \tsU {\widetilde {\sU}}
\define \th {\widetilde h}

\def\pd {\partial}

\def \Dx1 {\frac{\pd}{{\pd} x_1}}
\def \Dy1 {\frac{\pd}{{\pd} y_1}}
\def \Dz1 {\frac{\pd}{{\pd} z_1}}
\def \Dx2 {\frac{\pd}{{\pd} x_2}}
\def \Dy2 {\frac{\pd}{{\pd} y_2}}
\def \Dz2 {\frac{\pd}{{\pd} z_2}}

\def\q {\quad} 

\def\mapdiagr#1{\Big\searrow\rlap{$\raise 5pt\vbox{{\hbox{$\mkern -15mu\scriptstyle#1$}}}$}}   

\def\mapdiagl#1{\llap{$\raise 5pt\vbox{{\hbox{$\scriptstyle#1\mkern
-15mu$}}}$}\Big\swarrow}              

\def\Mapdiagr#1{\nearrow\rlap{$\lower 5pt\vbox{{\hbox{$\mkern
-15mu\scriptstyle#1$}}}$}} 

\def\Mapdiagl#1{\llap{$\lower 5pt\vbox{{\hbox{$\scriptstyle#1\mkern
-15mu$}}}$}\searrow} 

\def\Mapswr#1{\swarrow\rlap{$\lower 5pt\vbox{{\hbox{$\mkern
-15mu\scriptstyle#1$}}}$}}              

\def\Mapnwl#1{\nwarrow\rlap{$\lower 5pt\vbox{{\hbox{$\mkern
-15mu\scriptstyle#1$}}}$}}

\def \inj {\hookrightarrow}

\define \Rhook {\hookrightarrow}

\def \half {\raise1pt\hbox{$\scriptstyle
        \frac{1}{2}\displaystyle$}}

\def \x{{\sl X}\llap{$\mkern -2mu {\scriptstyle -}$}}
\def \Symm {\operatorname{Sym}}

\define \Kod {\operatorname{Kod}}
\define \dimension {\operatorname{dim}}
\define \codim {\operatorname{codim}}
\define \contr {\operatorname{contr}}
\define \rk {\operatorname{rank}}
\define \im {\operatorname{im}}
\define \Mor {\operatorname{Mor}}
\define \Cl {\operatorname{Cl}}
\define \Hilb {\operatorname{Hilb}}
\define \degree {\operatorname{deg}}
\define \mult {\operatorname{mult}}
\define \Aut {\operatorname{Aut}}
\define \NS {\operatorname{NS}}
\define \Gal {\operatorname{Gal}}
\define \ch {\operatorname{char}}
\define \Jac {\operatorname{Jac}}
\define \Km {\operatorname{Km}}
\define \Sec {\operatorname{Sec}}
\define \Stab {\operatorname{Stab}}
\define \Br {\operatorname{Br}}
\define \inv {\operatorname{inv}}
\define \tr {\operatorname{tr}}
\define \Frob {\operatorname{Frob}}
\define \Symn {\operatorname{Sym}^n}
\define \Ev {\sE^\vee}
\define \ordp {\operatorname{ord}_p}
\define \Supp {\operatorname{Supp}}
\define \Ann {\operatorname{Ann}}
\define \disc {\operatorname{disc}}
\define \Lie {\operatorname{Lie}}
\define \embdim {\operatorname{embdim}}

\def\hod#1#2#3#4{\ensuremath{\if#30 H^{#2}({#1},{\cal O}_{#1}) \else 
 H^{#2}(#1,\Omega^{#3}\if\relax{#4}\relax_{#1}\else _{#1/#4}\fi)\fi}}

\begin{document}
\title{Siegel modular forms and the gonality of curves}

\author{N. I. Shepherd-Barron}

\address{D.P.M.M.S.\\
University of Cambridge\\
Cambridge CB3 0WB U.K.}
\email{nisb@dpmms.cam.ac.uk}

\maketitle

\begin{section}{Introduction}
Denote by $M_g$ and $A_g$ the coarse moduli spaces of genus $g$ curves
and principally polarized abelian $g$-folds, respectively, and 
let a superscript $S$ denote their Satake compactifications.
The main result of [CSB] is that the intersection
$M_{g+m}^S\cap A_g^S$, taken inside $A_{g+m}^S$, is far
from transverse; it contains the $m$th order infinitesimal
neighbourhood of $M_g^S$ in $A_g^S$. So
$\cup_m (M_{g+m}^S\cap A_g^S)$ is the formal completion
of $A_g^S$ along $M_g^S$ and there are no stable Siegel modular
forms that vanish along every moduli space $M_g$. 
The proof depends upon
the construction by Fay [F] of certain, very special, 
degenerating families of curves
for which Yamada [Y] could subsequently
establish a formula for (a part of) the derivative
of the period matrix as a certain explicit tensor of rank 
one. For an arbitrary degeneration the derivative is a 
tensor of higher rank, usually maximal, and it is 
more difficult to make use of this; cf. the assertion on 
p. 1 of the erratum to [G-SM].  
Interpreting Fay's formula in terms of the projective geometry
of the canonical model of the singular fibre then gives the result.

Here we prove similar results for the loci $V_{g,n,tot}$
in $M_g$ of $n$-gonal curves of genus $g$ with a point of total ramification,
for any fixed $n\ge 3$, as follows. 

\begin{theorem} (= \ref{Main}) There is no stable Siegel
modular form that vanishes on every locus $V_{g,n,tot}$.
In particular, there is no stable Siegel modular form
that vanishes on every trigonal locus.
\noproof
\end{theorem}

This sharpens
[CSB], but depends upon it.
For hyperelliptic curves, however, Codogni has shown [C] that
the story becomes very different. He has found many
millions of stable modular forms that vanish on the hyperelliptic locus
in every genus, for example, the difference $\Theta_P-\Theta_Q$
of two theta series where $P,Q$ are positive, even and unimodular
quadratic forms of rank $32$ with no roots.

Recall that a curve $C$ is $n$-gonal if there is a map
$C\to\P^1$ of degree $n$. If $g=2n-2$ is even, then a general curve
of genus $g$ is $n$-gonal in finitely many ways; if $g>2n-2$
then the $n$-gonal curves form a proper subvariety
(the Hurwitz scheme) $V_{g,n}$ of $M_g$.
The $n$-gonal curves for which the given map
to $\P^1$ has a point of total ramification
form the subvariety $V_{g,n,tot}$ mentioned above.
Its closure in $M_g^S$ will be denoted by $V_{g,n,tot}^S$.

Compared to the arguments in [CSB], the proof here depends 
upon combining Fay's construction with those 
by Schiffer to get certain variations
of a curve where what is essentially the derivative
of the period matrix can be 
calculated explicitly. Controlling the construction
of these Fay--Schiffer variations (see below) is crucial
in controlling the derivative. 
\end{section}
\begin{section}{Variations}
Suppose that $C$ is a curve (= compact Riemann surface) of
genus $g$, that $a,b,c$ are distinct points of $C$ and that
$z_a,z_b,z_c$ are local co-ordinates on $C$ at $a,b,c$
respectively. There are various
well known kinds of variation that can be
constructed from these data, and we recall some of them now.

The first is a \emph{Fay variation} of $C$ centred at
$(a,z_a;b,z_b)$. This is a particular proper morphism $\sC\to\Delta$
from a smooth complex surface to a disc such that
the fibre over $0$ is the nodal curve $C/(a\sim b)$
and for every $t\ne 0$ the fibre $\sC_t=C_t$
is of genus $g+1$. It is constructed as follows
[F, p. 50].

Fix $\delta>0$ with $\delta<<1$.
Let $D_{\delta^2}$ be a disc of radius $\delta^2$
and complex co-ordinate $t$. In $C\times D_{\delta^2}$
consider two closed subsets, one defined by
the inequality $|z_a|\le |t|/\delta$ and the other
by the inequality $|z_b|\le |t|/\delta$.
Delete these closed subsets from $C\times D_{\delta^2}$
to get the complex manifold $\sC^0$. There are open
subsets $U_a$ and $U_b$ of $\sC^0$ defined by
the further inequalities 
$|z_a|<\delta$ and $|z_b|<\delta$, respectively.

Let $S$ be the
open part of the complex surface with co-ordinates $X,Y$
defined by the inequalities $|X|,|Y|<\delta$. There
is a morphism $S\to D_{\delta^2}$ given by $t=XY$.
Now map $U_a$ and $U_b$ to $S$ by the formulae
$$X=z_a,\ Y=t/z_a,$$
$$X=t/z_b,\ Y=z_b$$
and then glue $\sC^0$ to $S$ via these maps; 
by definition, the result is $\sC$, and $\sC$ is provided
with a proper morphism to $\Delta=D_{\delta^2}$.

Another kind is a
\emph{Schiffer variation} of $C$ centred at $(c,z_c)$.
This is a particular proper morphism $\sC\to\Delta$ where now
all fibres are smooth of genus $g$. 
It is also constructed via a glueing procedure.

Start with $C\times D_{\delta^2/4}$ and delete the closed subset
defined by the inequality 
$|z_c|\le {\sqrt{|t|}}$
to obtain the complex manifold $\sC^0$.
In $\sC^0$ there is the open subset $V_c$ defined by
$${\sqrt{|t|}}<|z_c| <\delta-{\sqrt{|t|}}.$$
The principle of the argument says that, as $z_c$ goes once
around the circle $R$ of radius $\delta-{\sqrt{|t|}}$ and
centre $0$, so $w=z_c+t/z_c$ has exactly one zero inside
$R$, so that the image of $R$ in the $w$-plane
is a simple closed curve $\G(t)$ around $0$, and varies
smoothly with $t$ for $0\le |t|<\delta^2/4$.

Say that $D(t)$ is the open neighbourhood of $0$
with boundary $\G(t)$. 
Then $\cup_{0\le |t|<\delta^2/4}D(t)$ is an open submanifold
$V$ of $\C\times D_{\delta^2/4}$. Map $V_c$ to $V$ via
$w=z_c+t/z_c$; this is unramified, since the branch locus
is given by $z_c^2+t=0$, and glueing $\sC^0$ to $V$
via the map $V_c\to V$
that has just been constructed gives the Schiffer variation
of which we speak.

If now $(a_1,...,a_n)$ are distinct points of $C$ and $z_j=z_{a_j}$
is a local co-ordinate at each, then we can
simultaneously construct a Fay variation centred
at $(a_{n-1},z_{n-1};a_{n},z_{n})$ and a Schiffer variation
centred at $(a_1,z_1;...;a_{n-2},z_{n-2})$.
This is a proper map $f:\sC^+\to \Delta^{n-1}$, where now 
$\Delta^{n-1}$ is an
$(n-1)$-dimensional complex polydisc with co-ordinates
$t_1,t_2,...,t_{n-1}$, the map $f$ is smooth over the locus
$t_{n-1}\ne 0$ and the fibres over $t_{n-1}=0$ are copies of the
nodal curve $C/(a_{n-1}\sim a_n)$.
We call it the Fay--Schiffer variation of $C$ centred at
$(a_1,z_1;...;a_{n},z_{n})$.

\begin{theorem}\label{2.1}
With respect to a suitable fixed homology basis
and a correspondingly normalized basis
$\omega=(\omega_1,...,\omega_g)$ of the abelian differentials on $C$,
the period matrix $T(t)$ of $C_t$ can be written in $2\times 2$
block form as
$$T(t)=
\left[
{\begin{array}{cc}
{\tau +\sum_1^{n-1}t_j\sigma_j} & 
{AJ(t)+t_{n-1}s}\\
{{}^t(AJ(t)+t_{n-1}s)} & {\frac{1}{2\pi i}(\log t_{n-1}+c_1+c_2t_{n-1})}
\end{array}}
\right] +O(t^2)
$$
where for $1\le j\le n-2$ the matrix $\sigma_j$ is of rank $1$ and is given by
$$(\sigma_j)_{pq}=2\pi i\left(\frac{\omega_p}{dz_j}(a_j).
\frac{\omega_q}{dz_j}(a_j)\right),$$
the matrix $\sigma_{n-1}$ is of rank $2$ and is given by
$$(\sigma_{n-1})_{pq}=2\pi i
\left(\frac{\omega_p}{dz_{n-1}}(a_{n-1})\frac{\omega_q}{dz_n}(a_n)+
\frac{\omega_q}{dz_{n-1}}(a_{n-1})\frac{\omega_p}{dz_n}(a_n)\right),$$
${}^tM$ is the transpose of the matrix $M$,
$AJ(t)=AJ_0(a_n-a_{n-1})+\sum_{j=1}^{n-2}t_jAJ_j$,
$AJ_0$ is the Abel--Jacobi map $AJ_0(y-x)=\int_x^y\omega$ on $C$,
each $AJ_j$ is a holomorphic function of the
parameters $a_i,z_i$ for $i=1,...,n-2$, 
$s, c_1, c_2$ are holomorphic 
functions of 
the parameters $a_j,z_j$ in the construction
and $c_1$ also depends on $t_1,...,t_{n-2}$.
\begin{proof} Consider the Schiffer variation of $C$ centred at 
$(a_1,z_1;...;a_{n-2},z_{n-2})$. This gives a genus
$g$ family $\Gamma\to \Delta^{n-2}$ where $\Delta^{n-2}$
is an $(n-2)$-dimensional polydisc with co-ordinates
$t_1,...,t_{n-2}$ and the period matrix of $\Gamma_t$ is
$$\left[\tau+\sum_1^{n-2}t_j\sigma_j\right] +O(t^2).$$
(This formula is due to Patt [P].)
By construction, this Schiffer variation is trivial
outside neighbourhoods of the points
$a_1,...,a_{n-2}$, and so the points
$a_{n-1},a_n$ and the local co-ordinates $z_{n-1}, z_n$
come along for the ride. So now we make a Fay variation
of $\Gamma\to \Delta^{n-2}$ centred at 
$(a_{n-1},z_{n-1};a_n,z_n)$ to get 
$\sC\to\Delta^{n-1}$. The period matrix $T(t)$ of the
curve $C_t$ of genus $g+1$ is then
$$\left[
\begin{array}{cc}
{\tau+\sum_1^{n-2}t_j\sigma_j+t_{n-1}\sigma_{n-1}} 
& {AJ_{\Gamma_t}(a_{n}(t)-a_{n-1}(t)) + t_{n-1}s}\\
{{}^t\left(AJ_{\Gamma_t}(a_n(t)-a_{n-1}(t)) + t_{n-1}s\right)}
&{\frac{1}{2\pi i}\left(\log t_{n-1}+c_1t_{n-1}+c_2\right)}
\end{array}
\right] + O(t^2)$$
where the matrix $\sigma_{n-1}$ is as
described in the statement of the theorem
(this is the correct form, due to Yamada [Y],
of Fay's original, but incorrect, formula),
$AJ_{\Gamma_t}$ is the Abel--Jacobi map
for the curve $\Gamma_t$ and each of the terms 
$AJ_{\Gamma_t}(a_n(t)-a_{n-1}(t))$, $s$, $c_1$ and $c_2$ 
is a holomorphic function
of $t_1,...,t_{n-2}$ and the parameters
$a_1,...,a_{n-2}$ and $z_1,...,z_{n-2}$.
However, for $t_1=\cdots=t_{n-2}=0$
the family $\sC\to\Delta^{n-1}$ is just
the usual Fay variation of $C$ centred at
$(a_{n-1},z_{n-1};a_n,z_n)$, and so
the Abel--Jacobi term $AJ_{\Gamma_t}(a_n(t)-a_{n-1}(t))$
is independent of the $a_j$ and the $z_j$; 
$AJ_{\Gamma_0}(a_n(0)-a_{n-1}(0))=AJ_C(a_n-a_{n-1})$
and so
$$AJ_{\Gamma_t}(a_n(t)-a_{n-1}(t))=AJ_C(a_n-a_{n-1})+
\sum_1^{n-2}t_jAJ_j+O(t^2).$$
\end{proof}
\end{theorem}

Now suppose that $h:C\to B$ is a morphism of Riemann surfaces
of degree $n$, that $e\in B$ is a point over which $h$ is unramified
and that $h^{-1}(e)=\{a_1,...,a_n\}$. For any local co-ordinate
$z_e$ on $B$ at $e$, define the local co-ordinate $z_j$ on $C$
at $a_j$ to be the pull-back of $z_e$ restricted to a neighbourhood
of $a_j$.

Take the corresponding Fay--Schiffer variation
$\sC^+\to\Delta^{n-1}$ of $C$ centred at
$(a_1,z_1;...;a_{n},z_{n})$, and let $\sC\to\Delta$
be the one-parameter family obtained by restricting
$\sC^+\to\Delta^{n-1}$ to the diagonal disc
$\Delta$ in $\Delta^{n-1}$ defined by
$t_1=\cdots=t_{n-1}=t$. Let $\sB\to\Delta$ be the
Schiffer variation of $B$ centred at $(e,z_e)$.

\begin{proposition}\label{2.2} There is a degree $n$ morphism
$H:\sC\to\sB$ relative to $\Delta$ that at $t=0$ is 
the morphism $C/(a_{n-1}\sim a_n)\to B$ induced by $h$.
\begin{proof} The Schiffer variation $\sB\to \Delta$
is constructed by deleting a disc and then glueing in 
a new disc with co-ordinate $w=z_e+t/z_e$; the variation
$\sC\to \Delta$ is constructed by the same formula
except where the points $a_{n-1},a_n$ are identified over
$t=0$. Here we have a complex surface $S$ with co-ordinates
$X,Y$ with $XY=t$, and the glueing was given by
$X=z_{n-1}, Y=t/z_{n-1}$ and $X=t/z_n, Y=z_n$.
So to construct $H:\sC\to\sB$ it is enough to give
the map from $S$ to the $w$-disc. This is achieved
by writing $w=X+Y$.
\end{proof}
\end{proposition}

Note that for all $t$, including $t=0$, the morphism
$H_t:C_t\to B_t$ coincides with $h$ outside a union
of small open sets. In particular, the ramification
data of $H_t$ coincides with those of $h$ away
from this union.

\begin{proposition}\label{closed}\label{2.3}
$V_{g,n,tot}\times M_1$ lies in the
closure of $V_{g+1,n,tot}$.
\begin{proof}
Suppose that the curve $C$ is a point in $V_{g,n,tot}$,
that $f:C\to\P^1$ is of degree $n$ and that $f$
is totally ramified at $P\in C$. Say $f(P)=e$,
so that $f^{-1}(e)=n[P]$. Suppose also
that the curve $E$ is a point in $M_1$. Fix $Q\in E$,
and regard $E$ as an elliptic curve with origin $Q$.
Then choose a primitive $n$-torsion point $R$ on $E$,
so that $n[Q]\sim n[R]$ and there is a rational
function $h:E\to\P^1$ such that $h^{-1}(0)=n[Q]$
and $h^{-1}(\infty)=n[R]$. We assume, as we may,
that $e\ne 0,\infty$.

We shall construct a variation similar (but not identical)
to that described on pp. 37--41 of [F], omitting
the topological details.
Choose local co-ordinates $z_e$ and $z_0$ on $\P^1$ at $e$
and $0$, respectively. Then there is a local co-ordinate
$w_P$ on $C$ at $P$ with $z_e=w_P^n$ and a local co-ordinate
$w_Q$ on $E$ at $Q$ with $z_0=w_Q^n$.
Use these to construct variations $\sC\to\Delta$ and
$\sB\to\Delta$, where $\sB$ is obtained by glueing
$\P^1\times\Delta$ and $\P^1\times\Delta$ to the surface
$S_n=(X_nY_n=t^n)$ by
$$X_n=z_e,\ Y_n=t^n/z_e,$$
$$X_n=t^n/z_0,\ Y_n=z_0$$
and $\sC$ is obtained by glueing $C\times\Delta$
and $E\times\Delta$ to the surface
$S_1=(X_1Y_1=t)$ by
$$X_1=w_P,\ Y_1 = t/w_P,$$
$$X_1=t/w_Q,\ Y_1=w_Q.$$
Via the morphism $S_1\to S_n$ given by $X_n=X_1^n,\ Y_n=Y_1^n$
there is a morphism $\pi:\sC\to\sB$ obtained by
glueing the morphisms
$f\times1_\Delta:C\times\Delta\to\P^1\times\Delta$
and $h\times1_\Delta:E\times\Delta\to\P^1\times\Delta$.
Moreover, since $h\times1_\Delta$ is totally ramified
along $\{R\}\times\Delta$ and the variation $\sC\to\Delta$
is trivial outside neighbourhoods of $P$ and of $Q$,
the morphism $\sC_t\to\sB_t$ is totally ramified somewhere.
Since $\sB_t\cong\P^1$, the result is proved.
\end{proof}
\end{proposition}
\end{section}

\begin{section}{Modular forms vanishing on $V_{g,n,tot}$}
We fix an integer $n$ with $3\le n\le g-1$. We are especially
interested in those values of $n$ for which a general curve of
genus $g$ possesses at most finitely many $g^1_n$'s, 
so that $n\le g/2 +1$.
Then if $C$ is a non-hyperelliptic curve possessing a 
pencil $\Pi$ that is a complete $g^1_n$,
the linear span $\langle D\rangle$
of each element $D$ of $\Pi$ is a copy of $\P^{n-2}$,
and as $D$ varies over $\Pi$ these copies sweep
out a rational scroll $\Sigma(\Pi)$ of dimension $n-1$ in $\P^{g-1}$.
For example, if $n=3$ then $\Sigma(\Pi)$ is a surface
(and is the intersection of the quadrics that contain $C$).

Suppose that $G=G_{g+1}$ is a Siegel modular form on
$A_{g+1}$ such that the restriction $G\vert_{M_{g+1}}$
of $G$ to $M_{g+1}$ has multiplicity at least $m$
along $V_{g+1,n,tot}$. 
That is, $G$ and all its partial derivatives $F$ of order $\le m-1$
with respect to the entries $T_{pq}$ of a period
matrix $T$ in $\mathfrak H_{g+1}$
in directions tangent to $M_{g+1}$ vanish along
$V_{g+1,n,tot}$. We can define the Siegel $\Phi$-operator
on the derivatives by
$$\Phi(F)(\tau)=\lim_{t\to i\infty}F\left(
\begin{array}{cc}
{\tau} & {0}\\
{0} & {t}
\end{array}
\right).
$$

\begin{lemma}\label{derivatives}
$\Phi(F)$ is a derivative of $\Phi(G)$
of order $\le m-1$ in directions tangent to $M_g$
and vanishes along $V_{g,n,tot}$.
\begin{proof}
By construction, $\Phi(F)$ can be computed by
restricting to $A_g\times A_1$, then restricting to
$A_g\times\{j\}$ for some $j\in A_1$, and finally
letting $j\to\infty$. Since the intersection of $M_{g+1}$
and $A_g\times A_1$ certainly contains $M_g\times M_1$,
the first part of the lemma is proved. The second
part then follows from Proposition \ref{closed}.
\end{proof}
\end{lemma}

\begin{theorem}\label{scrolls} Under these assumptions,
the restriction $\Phi(G)$ of $G$ to $M_g$
has multiplicity at least $m+1$ along $V_{g,n,tot}$.
\begin{proof}
We need to show that $\Phi(F)$ is singular along $V_{g,n,tot}$.
Now $F$ has a Fourier expansion
$$F(T)=\sum_{X\in S_{g+1}}a(X)\exp\pi i \tr(XT),$$
where $T$ is a point in Siegel space $\mathfrak H_{g+1}$
and $S_n$ is the lattice of positive semi-definite
$n\times n$ symmetric matrices over $\Z$ whose diagonal 
is even.

Take a curve $C$ in $V_{g,n,tot}$, and choose any reduced
divisor $D=\sum_1^n a_j$ in the specified $g^1_n$ on $C$.
Let $h:C\to B=\P^1$ be the morphism defined by this $g^1_n$
and say that $D=h^{-1}(e)$ and that $h$ is totally ramified at
$P$. We have, according to Proposition \ref{2.2}, 
a $1$-parameter Fay--Schiffer variation $\sC\to\Delta$
of $C$ centred at $(a_1,z_1;...;a_n,z_n)$
with a degree $n$ morphism to the Fay--Schiffer
variation $\sB\to\Delta$ of $B$ centred at $(e,z_e)$.
Since $B=\P^1$, the variation $\sB\to\Delta$ is trivial,
so that for $t\ne 0$ the curve $C_t$ lies in $V_{g+1,n}$.
Moreover, because the variation is constructed to be trivial
outside a neighbourhood of $D$, the curve
$C_t$ lies in $V_{g+1,n,tot}$.

Now the argument follows [CSB] closely.

Take $T=T(t)$ to be the period matrix of
$C_t$ as above. Note that since $t_1=\cdots=t_{n-1}=t$,
we can re-arrange $c_1$ and $c_2$ so that both
of them are independent of $t$, and are holomorphic
functions only of the parameters $(e,z_e)$. Then 
$$F_{g+1}(T) =
\sum_{X\in S_{g+1}}a(X)\exp\pi i\sum_{p,q=1}^{g+1}
x_{pq}T_{pq}$$
where $X=(x_{pq})$.
Our aim is to examine the coefficient of $t$ in the expansion
of this expression in powers of $t$, so calculate modulo $t^2$.
Since $\exp 2\pi i T_{g+1,g+1}\equiv t.\exp c_1.\exp (c_2t)$ modulo $t^2$, 
it follows that
$$(F_{g+1})(T) \equiv\sum_{x_{g+1,g+1}=0}{}+
\sum_{x_{g+1,g+1}=2}{}$$
modulo $t^2$,
since all terms with $x_{g+1,g+1}\ge 4$ vanish modulo $t^2$.
Here $\sum_{x_{g+1,g+1}=r}{}$ denotes the sum over
$X\in S_{g+1}$ with $x_{g+1,g+1}=r$,
for $r=0$ or $2$. 
Therefore, modulo $t^2$,
$$\sum_{x_{g+1,g+1}=0}{}\equiv\sum_{X\in S_g}a(X)
\exp\pi i\sum_{p,q=1}^gx_{pq}(\tau_{pq}+t\sigma_{pq})$$
and
\begin{eqnarray*}
\sum_{x_{g+1,g+1}=2}{}\equiv
& t.\exp c_1.\sum_{X\in S_{g+1},x_{g+1,g+1}=2}a(X)\\
&.\exp \left(2\pi i\sum_{p=1}^gx_{p,g+1}\int_{a_{n-1}}^{a_n}\omega_p\right)
.\exp \left(\pi i\sum_{p,q=1}^g x_{pq}\tau_{pq}\right).
\end{eqnarray*}
So the coefficient of $t$ is $A+B\exp c_1,$ where
$$A=\sum_{x_{g+1,g+1}=0}a(X)
\left(\pi i\sum_{p,q=1}^gx_{pq}\sigma_{pq}\right)
\left(\exp\pi i\sum_{p,q=1}^gx_{pq}\tau_{pq}\right),$$
$$B=\sum_{x_{g+1,g+1}=2}a(X)
\left(\exp 2\pi i\sum_{p=1}^gx_{p,g+1}
\int_{a_{n-1}}^{a_n}\omega_p\right)
\left(\exp \pi i\sum_{p,q=1}^g x_{pq}\tau_{pq}\right).$$
By assumption, $A+B\exp c_1$ vanishes identically.

Now rescale the local co-ordinate $z_e$; that is, given
any non-zero scalar $\lambda$, replace $z_e$ by
$\lambda^{-1}z_e$.
Such a rescaling will produce a different family $\sC\to \Delta$
with $C_t$ in $V_{g+1,n,tot}$ for all $t\ne 0$,
but the quantity $A+(\exp c_1)B$ will still vanish
for the rescaled family. Moreover,
$B$ is invariant under this
rescaling, as is revealed by a cursory inspection.
Also $c_1$ is a holomorphic function of $\lambda$
because the entries of a period matrix are holomorphic
functions of the parameters.

\begin{lemma} $A=B=0$.
\begin{proof}
From the description above of $\sigma_{pq}$,
this rescaling multiplies $\sigma_{pq}$
by $\lambda^2$, so that $A$ can be written as
$$A=C\lambda^2$$
with $C$ independent of $\lambda$. So we have an identity
$$C\lambda^2=-B\exp(c_1(\lambda))$$
of holomorphic functions on the $1$-dimensional
algebraic torus $\GG_m=\Sp\C[\lambda^\pm]$,
where $B,C$ are constant functions on $\GG_m$. 
The result follows at once.
\end{proof}
\end{lemma}

Now $A$ can also be written as
\begin{eqnarray*}
A=&\frac{\partial}{\partial t}\bigg\vert_{t=0}\left(
\sum_{X\in S_g}a(X)\exp\pi i\sum_{pq,=1}^g x_{pq}(\tau_{pq}+
t\sigma_{pq})\right)\\
&=\frac{\partial}{\partial t}\bigg\vert_{t=0}F_g(\tau+t\sigma).
\end{eqnarray*}
That is, $\sigma$ lies in the Zariski tangent space
$H$ at the point $\tau$ to the divisor in $\frak H_g$ defined
by the function $F_g$.
It is important to note that, from this description, $H$
depends upon $C$ but
is independent of any of the other parameters (points, local
co-ordinates) used to construct the variation.
Thus $H$ contains every $\sigma$ that arises from different
choices of these other parameters.

Assume that $C$ has no non-trivial automorphisms. Then there 
are the standard classical
natural identifications of tangent spaces to moduli given by
$$T_{[C]}M_g=H^0({\Omega^1_{C}}^{\otimes 2})^\vee,$$
$$T_{[C]}A_g =\Symm^2 H^0({\Omega^1_{C}})^\vee.$$
The inclusion $T_{[C]}M_g\inj T_{[C]}A_g$
is dual to the natural multiplication (which is surjective,
by Max Noether's theorem)
$\Symm^2 H^0({\Omega^1_{C}})\to H^0({\Omega^1_{C}}^{\otimes 2}).$
So the vector space of quadrics in $\P^{g-1}$ can be regarded
as the space of linear forms on $T_{[C]}A_g$, and then 
$T_{[C]}M_g$ is the subspace of $T_{[C]}A_g$
defined by the vanishing of those quadrics in $\P^{g-1}$
that contain $C$.

We know that the tangent space
$H$ to the divisor $(F_g=0)$ at the point $\tau$
in $\mathfrak H_g$ contains every matrix
$\sigma$ that arises as above.
Projectivize, and use the classical descriptions
above of the tangent spaces to moduli. Then (the projectivization 
of)
$H$ is a hyperplane in $\P(\Symm^2H^0(C,K_C))^\vee$
that contains every point $\sigma(n-1,n)=\sigma=(\sigma_{pq})$ of the form 
$$\sigma_{pq}=\left(\omega_p(a_n)\omega_q(a_{n-1})+
\omega_q(a_n)\omega_p(a_{n-1})\right)
+\sum_{j=1}^{n-2}\omega_p(a_j)\omega_q(a_j),$$
where we have omitted a factor of $2\pi i$ and
the factors of $dz_e$ that should appear
as denominators. We can also regard $H$ as a quadric
in the $\P^{g-1}$ in which $C$ is canonically embedded,
and then what we have to prove is that $H$ contains $C$. 

We shall in fact prove a stronger statement,
namely that $H$ contains the scroll $\Sigma(\Pi)$
(which certainly contains $C$) that is mentioned
in the first paragraph of this section.
 
In $\P^{g-1}$, any element $D=\sum_{j=1}^n a_j$ of the 
given pencil $\Pi$ spans a copy
$L=L_D$ of $\P^{n-2}$; the points $a_1,\ldots,a_n$ are,
therefore, in general position in $L$. Regard
$L$ as the projectivization
of an $(n-1)$-dimensional vector space $W$ and the points
$a_j$ as projectivizations of vectors 
$w_j=(\omega_1(a_j),\ldots,\omega_g(a_j))$ in $W$. Consider 
the second Veronese embedding $Ver_2(L_D)$ in a copy
$\P^N_D$ of $\P^N$, where
$N+1=n(n-1)/2$ and $\P^N_D$ is a linear subspace 
of the projectivized tangent space $\P(T_{[C]}A_g)$. Then
$H$ contains the point (in the projectivization
of $\Symm^2W$)
$$\sigma_{n-1,n}=\sigma_{n-1,n}(w_1,...,w_n)=w_{n-1}w_n+\sum_1^{n-2}w_j^2;$$
the same argument shows that $H$ also contains
every other point $\sigma_{k,l}$, for $k<l$,
that is obtained from $\sigma_{n-1,n}$ by permutation
of the vectors $w_1,...,w_n$. The $\sigma_{k,l}$ form a set of
$N+1$ points in $\P^N_D$.

\begin{lemma}\label{span} These $N+1$ points span $\P^N_D$.
\begin{proof} 
The vectors $w_j$ lie in a fixed vector space $\C^g$ and the
subset $\{w_1,...,w_n\}$ of $\C^g$
spans an $(n-1)$-dimensional subspace $W$ of $\C^g$.
Moreover, the vectors $w_j$,
and the subspace $W$ that they span, 
depend upon the choice of normalized
basis $(\omega_1,...,\omega_g)$ of $H^0(C,\Omega^1_C)$.
The normalized bases form a Zariski dense subset under the
$GL_g(\C)$-torsor that is the set of all bases of
$H^0(C,\Omega^1_C)$, so the collection of
subsets $\{w_1,...,w_n\}$ is Zariski dense 
in the symmetric product $(\C^g)^{(n)}$.
(Recall that $n<g$.)

Now suppose that the $\sigma_{kl}$ fail to span
$\Symm^2W$. That is, they are linearly
dependent in $\Symm^2(\C^g)$. Then, for every $n$ vectors
$w_1,...,w_n$ in $\C^g$ that are linearly dependent, the quantities
$\sigma_{kl}(w_1,...,w_n)$ are linearly dependent.
In particular, this is the case if
$\sum w_i=0$ and it is therefore enough
to prove Lemma \ref{span} under this additional hypothesis.

Then $W$ is the irreducible $(n-1)$-dimensional representation of 
the symmetric group $\mathfrak S_n$ as a Coxeter group of type 
$A_{n-1}$. Let $\mathbbm{1}$ denote the trivial $1$-dimensional 
representation, so that $W\oplus\mathbbm{1}$ is the standard 
permutation representation $V$ with standard basis
$(v_1,...,v_n)$ and $\mathbbm{1}$ is the line
generated by the first elementary symmetric function
$e_1=\sum v_i$. Let $\pi:V\to W$ be the projection,
so that $\pi(v_i)=w_i$. Note that in $\Symm^2V$
the kernel of the map induced by $\pi$ (which we still
denote by $\pi$) is just $e_1V$, where
$e_i=e_i(v_1,...,v_n)$ is the $i$th
elementary symmetric function.

Write $\tau_{kl}=\sigma(v_1,...,v_n)$,
so that $\pi(\tau_{kl}=\sigma_{kl}$.
Note that $\tau_{kl}=v_kv_l-v_k^2-v_l^2+p_2$,
where $p_2=\sum_j v_j^2 = e_1^2-2e_2.$
Let $T\subset \Symm^2V$ be the subspace spanned by the
$\tau_{kl}$.

\begin{lemma} The $\tau_{kl}$ are linearly independent.
\begin{proof}
Suppose that $\sum_{k<l}\lambda_{kl}\tau_{kl}=0$.
Then
$$\left(\sum_{p<q}\lambda_{pq}\right)p_2 = -\sum_{k<l}\lambda_{kl}
(v_kv_l-v_k^2-v_l^2).$$
On the RHS the coefficient of $v_kv_l$ is $-\lambda_{kl}$
and on the LHS it is $-2\sum_{p<q}\lambda_{pq}.$
So all the $\lambda_{kl}$ are equal, say to $\lambda$,
and then $2\binom{n}{2}\sum_{p<q}\lambda_{pq}=\lambda.$
So $\lambda=0$, as required.
\end{proof}
\end{lemma}

\begin{lemma} $T$ has zero intersection with $e_1V$.
\begin{proof}
Suppose that
$$0\ne \sum_{k<l}\lambda_{kl}\tau_{kl}\in e_1V.$$
That is,
$$\sum_{k<l}\lambda_{kl}
(v_kv_l-v_k^2-v_l^2+p_2)=e_1v=(v_1+\cdots +v_n)
(\alpha_1v_1+\cdots +\alpha_nv_n).$$
So, for $k<l$, we have $\lambda_{kl}=\alpha_k+\alpha_l$.
Define $\lambda_{kl}=\lambda_{lk}$ for $k>l$,
and set $\lambda_{kk}=0$. Then
$$LHS=\frac{1}{2}\sum_{k,l}\lambda_{kl}v_kv_l
-\frac{1}{2}\sum_{k,l}\lambda_{kl}v_k^2
-\frac{1}{2}\sum_{k,l}\lambda_{kl}v_l^2
+\frac{1}{2}\sum_{k,l}\lambda_{kl}p_2$$
and the coefficient of $v_k^2$ on the LHS is
$$-\frac{1}{2}\sum_{l}\lambda_{kl}
-\frac{1}{2}\sum_{l}\lambda_{kl}
+\frac{1}{2}\sum_{k,l}\lambda_{kl}$$
while on the RHS the coefficient is
$\alpha_k$.
Therefore
$$-\sum_q\lambda_{kq}+\frac{1}{2}\sum_{p,q}\lambda_{pq}=\alpha_k.$$
Replace $k$ by $l$ and add: the result is
$$-\sum_q\lambda_{kq}-\sum_r\lambda_{lr}+
\sum_{p,q}\lambda_{pq}=\lambda_{k,l}$$ for $k\ne l$.
Now fix $l$ and sum over all $k\ne l$ to get
$$-\sum_k\sum_q\lambda{kq}+\sum_q\lambda_{lq}
-(n-1) \sum_r\lambda_{lr}
+(n-1)\sum_{pq}\lambda_{pq}=\sum_k\lambda_{kl}.$$
Hence
$$(n-1)\sum_{pq}\lambda_{pq}-(n-1)\sum_r\lambda_{lr}
+\sum_k\lambda_{lk}=\sum_k\lambda_{lk},$$
so that $\sum_r\lambda_{lr}$ is independent of $l$.
Then $\lambda_{kl}$ is independent of $l$, and so of $k$,
so that
$$\sum_{k<l}(v_kv_l-v_k^2-v_l^2+p_2) =\alpha e_1^2.$$
Then
$e_2-\sum_{k\ne l}v_k^2+\binom{n}{2}p_2=\alpha e_1^2,$
which is an immediate contradiction.
\end{proof}
\end{lemma}

Now we can complete the proof of Lemma \ref{span}.
By the previous lemma, $T$ injects into $\Symm^2W$.
Since both have the same dimension, namely, $\binom{n}{2}$,
$T$ maps onto $\Symm^2W$, which is exactly
what was wanted.
\end{proof}
\end{lemma}

It follows that $H$ contains $\P^N_D$, and therefore contains 
$Ver_2(L_D)$ for every reduced divisor $D$ in $\Pi$, the $g^1_n$ under 
consideration.  So indeed $H$, when regarded as a quadric in 
$\P^{g-1}$, contains the rational scroll $\Sigma(\Pi)$.
\end{proof}
\end{theorem}

\begin{corollary}\label{n-gonal}
Assume that $n\ge 3$ and that $m\ge 1$. Then 
the intersection $V_{g+m,n,tot}^S\cap M_g$ contains
the $m$th order infinitesimal neighbourhood
of $V_{g,n,tot}$ in $M_g$.
\begin{proof}
Suppose that $\Phi$ is some modular form on $A_{g+1}$ such that
$(\Phi)_0\cap M_{g+1}$ is singular, with multiplicity
$m$, along $V_{g+1,n,tot}$. That is, $\Phi$ and all its derivatives
of order at most $m-1$, taken in directions along $M_{g+1}$,
vanish along $V_{g+1,n,tot}$. 

Suppose that $F$ is such a derivative. Then it follows from what 
we have shown that the restriction $F\vert_{M_g}$ is singular along 
$V_{g,n,tot}$. That is (and this is the content of Lemma \ref{derivatives}), 
the restriction 
$\Phi\vert_{A_g}$ of $\Phi$ to $A_g$
and all derivatives of $\Phi\vert_{A_g}$ of order at most $m$, taken
in directions along $M_g$, vanish along $V_{g,n,tot}$.
\end{proof}
\end{corollary}

\begin{theorem}\label{Main} Fix $n\ge 3$. Then there is no stable Siegel 
modular form that vanishes on the totally ramified
$n$-gonal locus $V_{g,n,tot}$ for every $g$.
\begin{proof} Suppose that $F$ is such a modular form.
Then, by Corollary \ref{n-gonal}, $F$ vanishes on $M_g$ for every $g$.
But the main result of [C-SB] is that then $F=0$.
\end{proof}
\end{theorem}

The main result of [G-SM] is that the
the Schottky form
$F=\Theta_{E_8^2}-\Theta_{D_{16}^+}$
(the difference of two theta series associated to the
positive even unimodular lattices $E_8^2$ and $D_{16}^+$ 
of rank $16$) that, by results of Schottky [S]
and Igusa [I1], [I2], defines $M_4$ inside $A_4$,
does not vanish along $M_5$. They prove further
that it cuts out the exactly trigonal locus $V_{5,3}$ in $M_5$, 
and does so with multiplicity $1$.

\begin{corollary} In genus $6$ the Schottky form $F$ does not
vanish along the totally ramified trigonal locus.
\begin{proof} Suppose that $F_6$ vanishes along $V_{6,3,tot}$.
Then, by Theorem \ref{scrolls}, 
the restriction $F_5\vert_{M_5}$ of $F_5$
to $M_5$ is singular along $V_{5,3,tot}$.
Then the trigonal locus $V_{5,3}$ is singular
along the subvariety $V_{5,3,tot}$. 
But the trigonal locus is smooth outside the hyperelliptic
locus, and we are done.
\end{proof}
\end{corollary}

For $g=6$ there is another 
subvariety of $M_g$ that is distinguished by 
the fact that the canonical model is not an intersection
of quadrics, namely the locus of plane quintics.
Our techniques, however, cannot let us decide
whether $F$ vanishes along this locus; more
generally, they cannot handle $g^r_d$'s with $r\ge 2$.

\end{section}
\begin{section}{The even genus case}
Suppose that $g=2(n-1)$ is even. Then a general
curve of genus $g$ has a finite, but non-zero, number
of $g^1_n$'s, while the locus $V_{g+1,n}$ is an
irreducible divisor in $M_{g+1}$ (and a general curve
in $V_{g+1,n}$ has a unique $g^1_n$).

Fix a general curve $C$ of genus $g=2(n-1)$, and let
$\Pi_1,...,\Pi_r$ be the $g^1_n$'s on it.
(The number $r$ is a known function of $g$, but all
we need is that $r\ge 4$ when $g\ge 6$.)
As above, the members of each $\Pi_i$ 
sweep out a scroll $\Sigma_i=\Sigma(\Pi_i)$ in
$\P^{g-1}$ that contains $C$.

\begin{lemma} 
If $g\ge 6$, then there is no quadric
in $\P^{g-1}$ that contains every $\Sigma_i$.
\begin{proof}
Choose any $a\in C$. For every $i$ there is a unique
$D_i\in\Pi_i$ passing through $a$. Say
$D_i=a+\sum_{j=2}^nb_{ij}$ and $L_i=\langle D_i\rangle$. 
Suppose that there is a hyperplane $H$ in $\P^{g-1}$
that contains each $L_i$; then
$$H.C\ge a +\sum_{ij}b_{ij},$$
so that $2g-2\ge 1 +r(n-1)$. Since
$r\ge 4$ this is impossible, and there is no
such hyperplane. Since the $L_i$ are linear,
this means that $\cup L_i$ has embedding dimension
$g-1$ at $a$. 

Now suppose that $Q$ is a quadric that contains
every $\Sigma_i$. Then $Q$ contains $\cup L_i$,
and so has embedding dimension $g-1$ at every point
of $C$. However, the singular locus
of a quadric is linear, and we are done.
\end{proof}
\end{lemma}

This is false for $g=4$; there are two $g^1_3$'s, but
the scrolls $\Sigma_1$ and $\Sigma_2$
coincide, and are the unique quadric containing $C$.

\begin{theorem} 
The $n$-gonal divisor $V_{g+1,n}$ in $M_{g+1}$ has 
contact with $A_g$ along $M_g$.
\begin{proof} 
We need to show that for any modular form $F=F_{g+1}$ on
$A_{g+1}$ that vanishes along $V_{g+1,n}$,
the restriction $F_g$ of $F$ to $A_g$ is singular
along $M_g$. But this follows from the proof
of Theorem \ref{scrolls} (the entire proof, except
for the final paragraph): if $F_g$ is smooth on $A_g$ at
the point $[C]$ of $M_g$, then the tangent hyperplane
$H$ to $A_g$ corresponds, if it is non-zero, 
to a quadric in $\P^{g-1}$
that contains every scroll $\Sigma(\Pi_i)$.
But we have just seen that there is no such quadric.
\end{proof}
\end{theorem}

\end{section}
\bigskip
I am very grateful to Giulio Codogni for many
conversations on these subjects and to the organizers of the Edinburgh 
meeting on modular forms in November 2012 for the stimulating
environment that led to the writing of this paper.
\bibliography{alggeom,ekedahl}

\providecommand{\bysame}{\leavevmode\hbox to3em{\hrulefill}\thinspace}
\begin{thebibliography}{EGAIII:2}







\bibitem[C]{C}
G.~Codogni, \emph{Non-perturbative Schottky problem and stable equations for the hyperelliptic locus},
arXiv:1306.1183

\bibitem[CSB]{CSB}
G.~Codogni and N.I.~Shepherd-Barron, \emph{The non-existence
of stable Schottky forms}, arXiv:1112.6137; Compositio Math., to appear.





\bibitem[F]{F}
J.~Fay, \emph{Theta functions on Riemann surfaces}, LNM \textbf{352},
Springer, 1973.

%


\bibitem[G-SM]{G-SM}
S.~Grushevsky and R.~Salvati~Manni, \emph{The superstring cosmological constant
and the Schottky form in genus $5$},
Am. J. Math. \textbf{133} (2011), 1007-1037 (erratum,
\textbf{134} (2012), 1139-1142), and arXiv:0809.1391v5.





\bibitem[I1]{I1}
J.-I.~Igusa, \emph{Schottky's invariant and quadratic forms}, in E. B. Christoffel Int.
Symp., Aachen (1981), 352-362.

\bibitem[I2]{I2}
\bysame, \emph{On the irreducibility of Schottky's divisor}, J. Fac. Sci. Univ. Tokyo Section 
IA Math., \textbf{28} (1981), 531-545.











\bibitem[S]{S}
F.~Schottky, \emph{Zur Theorie der Abelschen Functionen von vier Variabeln},
J. Reine Angew. Math. \textbf{102} (1888), 304-352.





\bibitem[Y]{Y}
A.~Yamada, \emph{Precise variational formulas for abelian differentials},
Kodai Math. J. \textbf{3} (1980), 114-143.

\end{thebibliography}
\bibliographystyle{pretex}
\end{document}